\documentclass[12pt,letterpaper,reqno]{amsart}
\usepackage{amscd}
\usepackage{latexsym}
\usepackage{amsfonts}
\usepackage{amssymb}
\newtheorem{theorem}{Theorem}[section]
\newtheorem{proposition}[theorem]{Proposition}
\newtheorem{lemma}[theorem]{Lemma}
\newtheorem{corollary}[theorem]{Corollary}
\theoremstyle{definition}
\newtheorem{definition}[theorem]{Definition}

\newcommand{\rk}{\mathrm{rk}}
\newcommand{\map}{\mathrm{map}}
\newcommand{\hocolim}{\mathrm{hocolim}}
\newcommand{\spaces}{\mathbf{Spaces}}
\newcommand{\Hom}{\mathrm{Hom}}
\newcommand{\Inn}{\mathrm{Inn}}

\newcommand{\res}{\mathrm{res}}
\newcommand{\Vect}{\mathrm{Vect}}
\newcommand{\Rep}{\mathrm{Rep}}
\newcommand{\im}{\mathrm{Im}}
\newcommand{\out}{\mathrm{Out}}
\newcommand{\Char}{\mathrm{Char}}
\newcommand{\Syl}{\mathrm{Syl}}
\begin{document}
\title[a quotient of the set ${[BG,BU(n)]}$]{A Quotient of the Set $[BG, BU(n)]$ for a Finite Group $G$ of Small Rank}
\author{Michael A. Jackson}
\address{The Ohio State University, 231 W. 18th Ave., Columbus, OH 43210, USA\\
mjackson@math.ohio-state.edu}

\maketitle

\begin{abstract}
Let $G_{p}$ be a Sylow 
$p$-subgroup of the finite group $G$ and let $\Char _{n}^{G}(G_{p})$ represent the 
set of degree $n$ complex characters of $G_{p}$ that are the 
restrictions of class functions on $G$. We construct a natural map $\psi _{G}
: [BG,BU(n)]\rightarrow \prod _{p| |G|} \Char _{n}^{G}(G_{p})$  and prove that $\psi _G$ is 
a surjection for all finite groups $G$ that do not contain a subgroup isomorphic to 
$(\mathbb{Z}/p)^3$ for any prime $p$. We show, furthermore, that $\psi _{G}$ is in fact a 
bijection for two types of finite groups $G$: those with periodic cohomology and those of odd 
order that do not contain a subgroup isomorphic to $(\mathbb{Z}/p)^3$ for any prime $p$. 

{\it 2000} MSC: 55R37, 55S35, 20J06, 20D15
\end{abstract}

\section{Introduction}
Our purpose is to investigate the relationship between the homotopy classes of maps 
from $BG$ for a finite group $G$ to $BU(n)$ and the degree $n$ characters of the 
Sylow $p$-subgroups of $G$. Throughout this paper we will let $G$ be a finite group,
$p$ a prime dividing the order of $G$, and $G_{p}$ a Sylow $p$-subgroup of $G$. 
In addition we will use $\Char _{n}(G_{p})$ to represent the set of 
degree $n$ complex characters of $G_{p}$ and will let $\Char _{n}^{G}(G_{p})$ represent the 
subset consisting of those degree $n$ complex characters of $G_{p}$ that are the 
restrictions of class functions on $G$.

Recall that for a group $G$, $BG$ is the classifying space of $G$. Also $[BG,BU(n)]$ is 
the homotopy classes of maps from $BG$ to $BU(n)$, which may also be thought of as 
$\pi _0 \Hom (BG,BU(n))$ (see \cite{dz:mbcs}). Next we will define a natural map  
$$\psi _{G}:[BG,BU(n)] \rightarrow \prod _{p| |G|} \Char _{n}^{G}(G_{p}),$$ which 
will be discussed throughout this paper. To define this map, let us examine the 
following diagram where $\Rep (G,U(n))=\Hom (G, U(n))/\Inn (U(n))$:
\[\begin{CD} 
[BG,BU(n)] @>{\cong }>> \prod_{p| |G|}[BG,BU(n)^{\wedge}_{p}] @>{\res}>> 
\prod_{p| |G|}[BG_{p},BU(n)^{\wedge}_{p}]\\
@V{\bar{\psi} _{G}}VV		@V{\bar{\phi} _{G}}VV		@A{\cong}AA\\
\prod _{p| |G|} \Char _{n}(G_{p})@= \prod _{p| |G|} \Char _{n}(G_{p}) 
@<{\cong}<< \prod_{p| |G|}\Rep (G_{p},U(n)).\\
\end{CD}\]
Notice that spaces in the center and right of the top row contain $BU(n)^{\wedge}_{p}$, which 
is the $p$-completion of the space $BU(n)$. (For more information on $p$-completion, see 
\cite[Chap. VI]{bk:hlcl}.) 
The bijection in the upper left is a result of work by  Jackowski, McClure, and Oliver 
\cite{jmo:htcs}, while the bijection on the far right is a result of work by Dwyer and 
Zabrodsky \cite{dz:mbcs} and the bijection in the lower right is a basic result in 
representation theory. The restriction map $res$ is induced by the inclusion of the Sylow p-subgroups 
$G_{p}$ into $G$. We now let maps $\bar{\psi} _{G}$ and $\bar{\phi} _{G}$ be the maps that 
make the diagram commute. 
The image of $\bar{\psi} _{G}$ and $\bar{\phi} _{G}$ both 
lie in the subset $\displaystyle{\prod _{p| |G|} \Char _{n}^{G}(G_{p})\subseteq \prod _{p| |G|} 
\Char _{n}(G_{p})}$ (see Theorem \ref{thm:char}). So we will let 
$\psi _{G}$ and $\phi _{G}$ be the maps $\bar{\psi} _{G}$ and $\bar{\phi} _{G}$ respectively, 
with the range restricted to $\displaystyle{\prod _{p| |G|}\Char _{n}^{G}(G_{p})}$. Now that we 
have defined the map $\psi _{G}$, we state the following three theorems, which 
express the main results of this paper:

\begin{theorem}\label{thm:rk1}
If $G$ is a finite group that does not contain a rank two elementary abelian 
subgroup, then the natural mapping 
$$\psi _{G}:[BG,BU(n)] \rightarrow \prod _{p| |G|} \Char _{n}^{G}(G_{p})$$ 
is a bijection.
\end{theorem}

\begin{theorem}\label{thm:odd}
Let $G$ be a finite group that does not contain a rank three elementary abelian subgroup. 
If $|G|$ is odd, then the natural mapping 
$$\psi _{G}:[BG,BU(n)] \rightarrow \prod _{p| |G|}\Char _{n}^{G}(G_{p})$$
is a bijection.
\end{theorem}

\begin{theorem}\label{thm:rk2}
If $G$ is a finite group that does not contain a rank three elementary abelian subgroup, 
then the natural mapping 
$$\psi _{G}:[BG,BU(n)] \rightarrow \prod _{p| |G|}\Char _{n}^{G}(G_{p})$$ is a surjection.
\end{theorem}

The results in this paper build on work by Mislin and Thomas \cite{mt:hsfg}, who a prove similar 
result to Theorem \ref{thm:rk1} where $U(n)$ is replaced by $SU(2)$. This work is also 
related to work by Jackowski and Oliver \cite{jo:vbcs}. They look at the Grothendieck group 
of $\Vect (BG)$ and show that it is isomorphic to $\prod _{p| |G|}\mathcal{R}(G_{p})^G$, where 
$\mathcal{R}(G_p)$ is the complex representation ring of $G_p$ restricted to thee elements that 
are stable under the action of $G$. The present 
work is intended to be a step toward classifying finite groups that act freely on a finite 
CW-complex that is homotopy equivalent to the product of two spheres. The classification 
of such groups began with Adem and Smith \cite{as:ospc,as:pcga}. Using Theorem 
\ref{thm:rk2} and this author's thesis \cite{j:vbobg}, a finite group can be shown to act 
freely on a finite CW-complex that is homotopy equivalent to the product of two spheres by 
demonstrating the appropriate element of $\prod _{p| |G|}\Char _{n}^{G}(G_{p})$. Such 
an element must be in the product of characters, not virtual characters, which correspond 
to the representation ring as used by Jackowski and Oliver. A more complete explanation 
of this application of the present work will be explored in a subsequent paper. 

\section{A homology decomposition and conjugation families}\label{sec:hdcf}
We start by showing that we can work with each prime dividing the order of $G$ 
separately. Notice that for the diagram from Section 1, the maps, excluding maps 
from $[BG,BU(n)]$, can each be separated into a product over all $p| |G|$, yielding a separate 
diagram for each $p| |G|$:
\[\begin{CD} 
[BG,BU(n)^{\wedge}_{p}] @>{\res_{p}}>> [BG_{p},BU(n)^{\wedge}_{p}]\\
@V{\bar{\phi }_{G,p}}VV		@A{\cong}AA\\
\Char _{n}(G_{p}) @<{\cong }<< \Rep (G_{p},U(n)).\\
\end{CD}\]
Notice that the image of the map $\bar{\phi }_{G,p}$ is always contained in the subset 
$\Char _{n}^{G}(G_{p})\subseteq \Char _{n}(G_{p})$. So we will let 
$\phi _{G,p}$ be the map $\bar{\phi }_{G,p}$ with the range restricted to $\Char _{n}^{G}(G_{p})$. 
We will show the image under certain hypotheses is all of $\Char _{n}^{G}(G_{p})$  by looking at the 
map $res_{p}$ and by introducing a homology decomposition of $BG$.

We will use the subgroup decomposition as our homology decomposition, which 
can be described in the following manner. First 
let $\mathcal{C}$ be a collection of $p$-subgroups of $G$ closed under conjugation.
Let $\mathcal{O_{C}}$ be the $\mathcal{C}$-orbit category, which is the category with objects 
$G/P$ for $P\in \mathcal{C}$ and with $G$-maps as the morphisms. Let $\mathcal{I}$ be the 
inclusion 
functor from $\mathcal{O_{C}}$ to the category of $G$-spaces. Composing $\mathcal{I}$ with 
the Borel construction $(-)_{hG}$ gives a functor $\alpha_{\mathcal{C}}~:~\mathcal{O_{C}}~
\rightarrow~\spaces$. Notice that $\alpha_{\mathcal{C}}(G/P)$ has the homotopy type of $BP$ 
for any $p$-subgroup $P\subseteq G$.
This functor naturally induces another functor 
$b_{\mathcal{C}}: \hocolim \textrm{ }\alpha  
_{\mathcal{C}} \rightarrow BG$; therefore, $b_{\mathcal{C}}$ gives a homology decomposition of 
$BG$ if and only if $\mathcal{C}$ is an ample collection of subgroups of $G$ (see \cite{d:hdcs}).

At this point an ample collection of subgroups of $G$ will be discussed, starting with three 
definitions.
\begin{definition}
Let $P \subseteq G$ be a $p$-subgroup. $P$ is said to be \emph{p-radical} (or 
$p$-stubborn) if $N_{G}(P)/P$ has no non-trivial normal $p$-subgroups.
\end{definition}
\begin{definition}
Let $P \subseteq G$ be a $p$-subgroup. $P$ is said to be \emph{p-centric} 
if $Z(P)$ is a Sylow $p$-subgroup of $C_{G}(P)$.
\end{definition}
\begin{definition}[see \cite{g:hlsc}]
Let $P \subseteq G$ be a $p$-centric subgroup. $P$ is said to be \emph{principal p-radical}
 if $N_{G}(P)/P C_{G}(P)$ has no non-trivial normal $p$-subgroups.
\end{definition}

Now let the collection $\mathcal{C}$ be the set of all principal $p$-radical subgroups of $G$. 
Grodal has shown that this is an ample collection \cite{g:hlsc}, allowing for 
a homology decomposition of $G$. Also it should be noted that any Sylow p-subgroup of $G$ 
is contained in the collection $\mathcal{C}$ and that any principal $p$-radical subgroup of 
$G$ is necessarily $p$-radical.\\
\\
Next notice that the map $res_{p}$ factors through an inverse limit constructed via
the homology decomposition described above:
 \[ res_{p} : [BG, BU(n)_{p}^{\wedge}] \stackrel {\alpha _{p}}{\rightarrow} 
{\lim _{\leftarrow _{G/P \in \mathcal{O_{C}}}}} [BP,BU(n)_{p}^{\wedge}]
\stackrel {\beta _{p}}{\rightarrow} [BG_{p},BU(n)_{p}^{\wedge}]. \]
The map $\alpha _{p}$ is induced by restriction and the map $\beta _{p}$ is a projection onto 
a particular element since $G_{p}\in \mathcal{C}$. 
The following diagram then commutes where 
the maps $\pi ^{1}_{p}$ and $\pi ^{2}_{p}$ are projection onto the the set of degree $n$ 
characters of $G_{p}$ and onto the set of representations of $G_{p}$ respectively:

$$\begin{CD}
[BG, BU(n)_{p}^{\wedge}]@= [BG, BU(n)_{p}^{\wedge}]\\
@V{\alpha _{p}} VV   @V{\res_{p}} VV                                  \\
{\displaystyle {\lim _{\leftarrow _{G/P \in \mathcal{O_{C}}}}} [BP,BU(n)_{p}^{\wedge}]}
@>{\beta _{p}}>> [BG_{p}, BU(n)_{p}^{\wedge}]  \\
@A{\cong }AA		@A {\cong }AA \\
{\displaystyle {\lim _{\leftarrow _{G/P \in \mathcal{O_{C}}}}} \Rep (P,U(n))}@>{\pi ^{2}_{p}}>>
\Rep (G_{p},U(n))  \\
@V{\cong }VV   @V{\cong }VV                  \\
{\displaystyle{\lim _{\leftarrow _{G/P \in \mathcal{O_{C}}}}} \Char _{n}(P)} @>{\pi ^{1}_{p}}>>
\Char _{n}(G_{p}).\\
\end{CD}$$
\vspace{0.5 cm}\\
The bijections on the left side follow easily from the earlier discussion of $p$-groups. We 
notice that the composite of the entire right hand column is the map $\bar{\phi }_{G,p}$ 
whose image we want to find. From the diagram it is obvious that we can instead find the 
image of the map $\pi ^{1}_{p}$.\\
\\
In order to look closer at the image of the map $\pi ^{1}_{p}$, Alperin's fusion theory will be 
discussed. In particular the definition of weak conjugation family and an example will be given. 

\begin{definition}[Alperin \cite{a:siaf}]
Let $G$ be a finite group, $p$ a prime dividing $|G|$, and $G_{p}$ a Sylow $p$-subgroup of $G$. 
A set $\mathcal{F}$ of pairs $\{(H,T)\}$, where $H\subseteq G_{p}$ and $T\subseteq N_{G}(H)$, 
is called a \emph{weak conjugation family} provided that whenever $A$ and $B$ are 
subsets of $G_{p}$ and 
$B=g^{-1}Ag$ for $g\in G$, there are elements $(H_{1}, T_{1}),(H_{2},T_{2}),\ldots ,
(H_{n},T_{n})$ of $\mathcal{F}$ and elements $x_{1},x_{2},\ldots ,x_{n},y $ of $G$ such that
\begin{enumerate}
\item $B=(x_{1}x_{2}\cdots x_{n}y)^{-1}A (x_{1}x_{2}\cdots x_{n}y)$,
\item $x_{i}\in T_{i}$ for $1\leq i \leq n$ and $y\in N_{G}(G_{p})$, and
\item $A\subseteq H_{1}$, $(x_{1}x_{2}\cdots x_{i})^{-1}A (x_{1}x_{2}\cdots x_{i})
\subseteq H_{i+1}$ for $1\leq i \leq n-1$.
\end{enumerate}
\end{definition}
The example we will be using of a weak conjugation family is given in the following theorem 
of Goldschmidt. In order to state the theorem we need to give two more definitions.
\begin{definition}[Alperin \cite{a:siaf}]
Let $Q$ and $R$ be Sylow $p$-subgroups of a finite group $G$ and let 
$H=Q\cap R$. $H$ is said to be a \emph{tame intersection} if $N_{Q}(H)$ and $N_{P}(H)$ are 
Sylow $p$-subgroups of $N_{G}(H)$.
\end{definition}

\begin{definition}
A finite group $G$ is called \emph{$p$-isolated} if it contains a proper subgroup $H\subset G$ 
such that if $p||H|$ and for any $g\in G\setminus H$, $p\not ||H\cap gHg^{-1}|$. In this case 
$H$ is called a strongly $p$-embedded subgroup of $G$.
\end{definition}

\begin{theorem}[Goldschmidt {\cite[Theorem 3.4]{g:acffg}}]
Let $G$ be a finite group, $p$ a prime dividing $|G|$, and $G_{p}$ a Sylow $p$-subgroup of $G$. 
Let $\mathcal{F}_{0}$ be the 
set of all pairs $(H,N_{G}(H))$ where $H\subseteq G_{p}$ such that there exists a Sylow 
$p$-subgroup $P$ of $G$ with the following properties:
\begin{enumerate}
\item $H=G_{p}\cap P$ a tame intersection, 
\item $C_{G_{p}}(H) \subseteq H$, 
\item $H$ a Sylow $p$-subgroup of $O_{p',p}(N_{G}(H))$, and 
\item $H=P$ or $N_{G}(H)/H$ is $p$-isolated. 
\end{enumerate}
$\mathcal{F}_{0}$ is a weak conjugation family.
\end{theorem}

The next proposition will begin to relate Alperin's fusion theory with the inverse limit that 
has been used in this paper. This will allow us to find the image of the map $\pi ^{1}_{p}$.
\begin{proposition}
Let $G$ be a finite group, $p$ a prime dividing $|G|$, and $G_{p}$ a Sylow $p$-subgroup of $G$. 
The set $\mathcal{F}$ consisting of all 
pairs $(H,N_{G}(H))$, where $H\subseteq G_{p}$ is a principal $p$-radical subgroup of $G$, 
is a weak conjugation family.
\end{proposition}

\begin{proof} We will show that $\mathcal{F}$ is a weak conjugation family by showing that as a 
set it contains $\mathcal{F}_{0}$, the weak conjugation family of Goldschmidt. Let 
$(H,N_{G}(H))\in \mathcal{F}_{0}$. We first notice that $H$ is a $p$-centric subgroup 
of $G$ since $C_{G_{p}}(H)$ is a Sylow $p$-subgroup of $C_{G}(H)$ by the fact that $H$ is 
a tame intersection. Since $H$ is a Sylow $p$-subgroup of $O_{p',p}(N_{G}(H))$ and $H$ 
is $p$-centric, $H$ must also be principal $p$-radical (see \cite[Remark 10.12]{g:hlsc}). \end{proof}

\begin{theorem}\label{thm:char}
Let $G$ be a finite group, $p$ a prime dividing $|G|$, and $G_{p}$ a Sylow $p$-subgroup of G. 
If $\mathcal{C}$ is the collection of all $p$-subgroups $H\subseteq G$ that are prinicipal 
$p$-radical, then the projection map
$$ \pi ^{1}_{p}:{\lim _{\leftarrow _{G/P \in \mathcal{O_{C}}}}} \Char _{n}(P) \rightarrow 
\Char _{n}(G_{p})$$ is one to one and is onto the subset 
$\Char _{n}^{G}(G_{p})\subseteq \Char _{n}(G_{p})$.
\end{theorem}

\begin{proof} First we will show that the image of $\pi ^{1}_{p}$ is contained in $\Char _{n}^{G}(G_{p})$. 
Suppose that 
$\gamma  \in  {\displaystyle \lim _{\leftarrow_{G/P \in \mathcal{O_{C}}}}}\Char _{n}(P)$. 
Define a map $\pi _{P}$ for each $P\in \mathcal{C}$ as the projection from the inverse limit 
of the complex character of $P$. Notice that $\pi ^{1}_{p}=\pi _{G_{p}}$. Next define 
$\chi _{P}= \pi _{P} (\gamma )$ and let $\chi $ mean $\chi _{G_{p}}$. Let $\mathcal{F}_{p}$ 
be the set consisting of all pairs $(H,N_{G}(H))$ where $H\in \mathcal{C}$ and 
$H\subseteq G_{p}$. By the last proposition, $\mathcal{F}_{p}$ is a weak conjugation family.

Suppose that $a,b \in G_{p}$ such that 
$\exists \; g\in G$ with $b=gag^{-1}$. By the definition of weak conjugation family, 
there exists $H_{1}, H_{2}, \ldots , H_{m}\in 
\mathcal{C}$ with $H_{i}\subseteq G_{p}$ for $1\leq i \leq n$, and $x_{1},x_{2},\ldots x_{m},y\in G$ 
such that 
\begin{itemize}
\item $b=(x_{1}x_{2}\cdots x_{n}y)^{-1}a (x_{1}x_{2}\cdots x_{n}y)$, 
\item $x_{i}\in N_{G}(H_{i})$ for $1\leq i \leq n$,
\item $y\in N_{G}(G_{p})$, and 
\item $a\subseteq H_{1}$, $(x_{1}x_{2}\cdots x_{i})^{-1}a(x_{1}x_{2}\cdots x_{i})\subseteq 
H_{i+1}$ for $1\leq i \leq n-1$. 
\end{itemize}
For each $1\leq i \leq n$ we notice 
that $\chi |_{H_{i}}$ respects fusion in $N_{G}(H_{i})$ and $\chi $ respects fusion in 
$N_{G}(G_{p})$. From this we see that the following three statements hold:
\begin{itemize}
\item $\chi (a)=\chi (x_{1}^{-1}ax_{1})$,
\item $\chi ((x_{1}x_{2}\cdots x_{i})^{-1}a (x_{1}x_{2}\cdots x_{i})) 
= \chi ((x_{1}x_{2}\cdots x_{i+1})^{-1}a (x_{1}x_{2}\cdots x_{i+1}))$
for each $1\leq i \leq n-1$, and
\item  $\chi ((x_{1}x_{2}\cdots x_{n})^{-1}a (x_{1}x_{2}\cdots x_{n}))=\chi (b)$.
\end{itemize}
Putting these statements together, we see that $\chi (a)=\chi (b)$; therefore, $\chi $ respects fusion 
in $G$ and so is contained in $\Char _{n}^{G}(G_{p})$.

To show that $\pi ^{1}_{p}$  is a one to one correspondence, it is enough to show the existence of 
an inverse mapping
$$\sigma  _{G_{p}}: \Char _{n}^{G}(G_{p}) \rightarrow {\lim _{\leftarrow_{G/P \in 
\mathcal{O_{C}}}}\Char _{n}(P}).$$
Let $\chi \in \Char _{n}^{G}(G_{p})$. We will define $\sigma _{G_{p}}$ by giving the characters 
$\chi _{P}$, which will be $\pi _{P}\circ \sigma_{G_{p}}(\chi )$ for each $P\in \mathcal{C}$. 
Fix $P\in \mathcal{C}$. There exists $g\in G$ such that $P\subseteq g^{-1}G_{p}g$. 
So let $\chi _{P}(h)=\chi (ghg^{-1})$. Since $\chi $ respects fusion in $G$, it is easy to see 
that this definition of $\chi_{P}=\pi _{P}\circ \sigma_{G_{p}}(\chi )$  is well defined. 
This definition then gives the following definition of $\sigma  _{G_{p}}$:
$$\sigma  _{G_{p}}( \chi ) = {\lim _{\leftarrow_{G/P \in \mathcal{O_{C}}}}}\chi _{P} .$$

It is obvious from the definition that $\pi ^{1}_{p}\circ \sigma _{G_{p}}$ is the identity mapping. 
We now have to show only that $\sigma _{G_{p}} \circ \pi ^{1}_{p}$ is the identity mapping. 
Let $\gamma $ be an element of the inverse limit and let $\chi =\pi _{G_{p}}(\gamma )$. Given 
$P\in \mathcal{C}$, observe that $\pi _{P}(\gamma )$ must be $\chi _{P}$ defined above by the 
nature of the $\mathcal{C}$-orbit category $\mathcal{O_{C}}$. This observation shows that 
$\pi _{P} \circ \sigma _{G_{p}} \circ \pi ^{1}_{p}=\pi _{P}$ for each $P\in \mathcal{C}$; therefore, 
$\sigma _{G_{p}} \circ \pi ^{1}_{p}$ must be the identity mapping. \end{proof}

Applying Theorem \ref{thm:char} shows that the map $\pi ^{1}_{p}$ is a bijection when the 
range is restricted to $\Char _{n}^{G}(G_{p})$. Using this result we see that we can use the maps 
$\phi _{G,p}$, $\phi _{G}$, and $\psi _{G}$ instead of the corresponding maps $\bar{\phi }_{G,p}$, 
$\bar{\phi }_{G}$, and $\bar{\psi }_{G}$. We also see that if the map $\alpha _{p}$ is an 
injection or a surjection, so is the map $\phi _{G,p}$. In order to study the map 
$\alpha _{p}$, we must apply obstruction theory.

\section{Obstruction theory}

Recall the functor $\alpha _{\mathcal{C}}:\mathcal{O_{C}}\rightarrow\spaces $, 
which we encountered in Section~\ref{sec:hdcf}. Fixing an element 
$$\gamma  = (\gamma _{(G/P)} )_{G/P \in \mathcal{O_{C}}} \in 
{\lim _{\leftarrow_{G/P \in \mathcal{O_{C}}}}}[\alpha _\mathcal{C} (G/P),BU(n)^{\wedge }_{p}],$$ 
we define functors $\mathcal{D}_n : \mathcal{O_C} \rightarrow \mathcal{A}b$ for each $n\geq 1$ 
by letting $$\mathcal{D}_n (G/P)=\pi _n (\map (\alpha _{\mathcal{C}} (G/P) , BU(n) _p ^{\wedge })
_ {\gamma _{(G/P)} }).$$
At this point in our discussion we will use a theorem by Jackowski, McClure, and 
Oliver \cite{jmo:htcs}. This particular theorem explains where 
the obstructions will lie in our examination of the map $\alpha _p$.

\begin{theorem}[Jackowski, McClure, and Oliver \cite{jmo:htcs}]\label{thm:jmo}
Fix an element \newline $\gamma  \in {\lim _{\leftarrow_{G/P \in \mathcal{O_{C}}}}}[\alpha 
_\mathcal{C} (G/P),BU(n)^{\wedge }_{p}]$.
$\gamma \in \im (\alpha _p ) $ if the groups $\lim ^{n+1} (\mathcal{D}_n)$ vanish for all 
$n\geq 1$ and $\alpha _p ^{-1}(\gamma ) $ contains at most one element if the groups 
$\lim ^n (\mathcal{D}_n)$ vanish for all $n \geq 1$. 
\end{theorem}

Since we will are trying to show that $\alpha _p$ is an injection or a surjection, we will use the following 
corollary to Theorem \ref{thm:jmo}:
\newpage
\begin{corollary}\label{cor:jmo}
If for each element $\gamma \in {\lim _{\leftarrow_{G/P \in \mathcal{O_{C}}}}}
[\alpha _\mathcal{C} (G/P),BU(n)^{\wedge }_{p}]$,\newline $\lim ^{n+1} (\mathcal{D}_n)$ vanish for all $n\geq 1$, then 
$\alpha _p$ is a surjection. On the other hand, if for each $\gamma$, $\lim ^n (\mathcal{D}_n)$ vanish for all 
$n \geq 1$, then $\alpha _p$ is an injection.
\end{corollary}

For additional discussion of the obstruction theory see \cite{b:hsso,bk:hlcl,g:hlsc,jmo:htcs,w:omfz}. 
In light of Corollary \ref{cor:jmo}, we are left to show that under the correct hypotheses the 
groups $\lim ^{n+1} (\mathcal{D}_n)$ and $\lim ^{n} (\mathcal{D}_n)$ vanish for 
$n\geq 1$. We will notice first that in the case where $n=1$, these groups vanish without 
any additional hypotheses. 
Recall that for any representation $\rho :P \rightarrow U(m)$ with $P$ a $p$-group, 
$$\pi _{n} (\mathcal{D_{C}}(G/P),B\rho  )\cong \pi _{n} (BC_{U(m)}(Im(\rho))) \otimes 
\mathbb{Z}^{\wedge}_{p} \textrm{ (see \cite{dz:mbcs}).}$$ Since the centralizer of a 
finite $p$-subgroup 
of $U(m)$ is the product of various $U(i)$, 
it is clear that $BC_{U(m)}(Im(\rho))$ is simply connected; therefore, for any element 
$\gamma $,
$\lim ^{2} (\mathcal{D}_1)=0=\lim^{1}  (\mathcal{D}_1)$. In the remaining sections we will 
discuss situations where the remaining higher limits are indeed trivial.

\section{Vanishing obstruction groups}
Before we proceed in showing when these higher limits vanish, we first need a couple of 
definitions. These will allow us to use standard theory first introduced by Jackowski, McClure, 
and Oliver \cite{jmo:hcsm,jmo:htcs}.
\begin{definition}
A functor $F:\mathcal{O_{C}} \rightarrow \mathit{Ab}$ is called an \emph{atomic} functor if 
it vanishes on all but possibly a single isomorphism class of objects. An atomic functor 
will be said to be concentrated on an isomorphism class on which it does not vanish.
\end{definition}
\begin{definition}[Definition 4.7 of \cite{jmo:htcs}]
For any prime $p$, finite group $G$, and 
$\mathbb{Z}_{(p)}(G)$-module $M$, let $F_{M}:\mathcal{O}_{\mathcal{S}_{p}(G)} \rightarrow 
\mathbb{Z}_{(p)}-mod$ be the atomic functor concentrated on the free orbit $G/1$ with 
$F_{M}(G/1)=M$. Set $\Lambda ^{\ast}(G;M) = {\lim}^{\ast } F_{M}$.
\end{definition}
The following lemma follows directly from work by Jackowski, McClure, and Oliver, in 
particular from the first proof of Lemma 5.4 in \cite{jmo:hcsm}.
\begin{lemma}
Suppose $\mathcal{C}$ is a subcollection of $\mathcal{S}_{p}(G)$ for a finite group $G$ 
such that for any functor $F:\mathcal{O}_{\mathcal{S}_{p}(G)} \rightarrow \mathbb{Z}_{(p)}-mod$, 
${\lim}^{\ast } F \cong {\lim}^{\ast } F|_{\mathcal{O_{C}}}$. 
Also assume that $F':\mathcal{O_{C}} \rightarrow \mathbb{Z}_{(p)}-mod$ is an atomic functor 
concentrated on the isomorphism class of $G/P$ for some $P\in \mathcal{C}$. Then 
$${\lim}^{\ast } F' \cong \Lambda ^{\ast }(N_{G}(P)/P;F'(G/P)).$$
\end{lemma}
We will be filtering our functor by atomic functors in order to see where 
$\Lambda ^{\ast }(\Gamma ,M)$ vanishes. The relevant lemma is a result of work by Grodal 
{\cite{g:hlsc}}. (See also {\cite[Prop. 4.11]{jmo:htcs}} and 
{{\cite[Proposition 5.8]{blo:hepcs}}.)
\begin{lemma}[See \cite{g:hlsc}]\label{lem:van}
Let $\Gamma $ be a finite group and $M$ a finitely generated 
$\mathbb{Z}_{(p)}(\Gamma )$-module. If 
$m$ is a non-negative integer such that $\rk _{p} (\Gamma )\leq m$, 
then $\Lambda ^{j}(\Gamma ,M)=0$ for any $j>m$. In particular, if 
$p^{m+1}\not ||\Gamma |$, 
then $\Lambda ^{j}(\Gamma ,M)=0$ for any $j>m$.
\end{lemma}
Lemma \ref{lem:van} gives us the following proposition when a functor is filtered by atomic 
functors.
\begin{proposition}
Let $G$ be a finite group and $p$ a prime number dividing $|G|$. Let $\mathcal{C}$ be 
the set of all subgroups of $G$ that are principal $p$-radical.
If $m$ is a positive integer such that for each $P\in \mathcal{C}$ 
$\rk _{p} (N_{G}(P)/P )\leq m$, then for any functor 
$F:\mathcal{O_{C}} \rightarrow \mathbb{Z}_{(p)}-mod$, ${\lim}^{j} F=0$ for any $j>m$.
\end{proposition}
Using this proposition in conjunction with the work of the previous section, 
we get the following result.
\begin{theorem}\label{thm:rnk}
Let $G$ be a finite group and $p$ a prime number dividing $|G|$. Let $\mathcal{C}$ be 
the set of all subgroups of $G$ that are principal $p$-radical.
If for each $P\in \mathcal{C}$ $\rk _{p} (N_{G}(P)/P )\leq 2$, then the map $\alpha _{p}$ 
described earlier is a surjection. If for each $P\in \mathcal{C}$ $\rk _{p} (N_{G}(P)/P )\leq 1$, 
then the map $\alpha _{p}$ is a bijection.
\end{theorem}
In the rest of this paper, we will use Theorem \ref{thm:rnk} and various results about finite groups 
to prove the three theorems presented in the introduction.

\section{Finite groups, part I}
We will start by proving Theorem \ref{thm:rk1}. We begin this process with the following 
lemma, which shows that the obstructions vanish when a Sylow 
$p$-subgroup of $G$ has a center of small index in itself. 
\begin{lemma}\label{lem:cent}
Let $G$ be a finite group, $p$ a prime number dividing $|G|$, and $G_{p}$ a Sylow 
$p$-subgroup of $G$. Let $m$ be an integer such that $p^{m}=[G_{p}:Z(G_{p})]$. If $P$ is a 
$p$-centric subgroup of $G$, then $p^{m}\not |[N_{G}(P):P]$. In particular, if $G_{p}$ is 
abelian, then the only $p$-centric subgroups of $G$ are the Sylow $p$-subgroups.
\end{lemma}
\begin{proof} Let $P\subseteq G$ be a $p$-centric subgroup, which is not a Sylow $p$-subgroup of 
$G$. We may assume $P\subset G_{p}$, which implies that $Z(G_{p})\subseteq C_{G}(P)
\cap G_{p}$. Let $H=~C_{G}(P)\cap G_{p}$, which is a Sylow $p$-subgroup of $C_{G}(P)$. 
It is clear that $Z(P)\subseteq H$, which implies that $Z(P)=H$ since $P$ is $p$-centric; therefore, 
$Z(G_{p}) \subseteq Z(P)\subseteq P$. Suppose that $Z(G_{p})=P$, then 
$G_{p}\subseteq C_{G}(P)$, contradicting $P$ being $p$-centric. So we see that $Z(G_{p})$ 
is strictly contained in $P$. Now by hypothesis, $p^{m} = [G_{p}:Z(G_{p})]$; thus, by the 
strict containment $Z(G_{p})\subset P$, $[G_{p}:P]< p^{m}$. This clearly implies 
$p^{m} \not | [G:P]$. Then since $[N_{G}(P):P]$ divides 
$[G:P]$, $p^{m}\not | [N_{G}(P):P]$. \end{proof}

At this point we will prove Theorem \ref{thm:rk1}:
\begin{proof} As was discussed earlier, it is enough to show that for each prime $p$ that divides $|G|$, 
$\phi _{G,p}$ is a bijection. 
We assume that $G$ does not contain a rank two elementary abelian subgroup. 
Let $G_{p}$ be a Sylow $p$-subgroup of $G$. If $G_{p}$ is cyclic, it is abelian; therefore, 
the only $p$-centric subgroups of $G$ are the Sylow $p$-subgroups. By 
Theorem \ref{thm:rnk}, the map $\alpha _{p}$ is a bijection as is the map 
$\phi _{G,p}$. On the other hand, if $G_{p}$ is not a cyclic group, then $p=2$ 
and $G_{p}$ is a generalized quaternion group. Suppose in this case
that $P\subseteq G_{p}$ is a $p$-centric subgroup of $G$. We see that either $P$ is the 
cyclic group of index two in $G_{p}$, $P$ is cyclic of order four and intersects the cyclic group 
of index two in $G_{p}$ in a cyclic group of order two, or $P$ itself is a generalized quaternion 
group. If $P$ is of index two in $G_{p}$, 
$4\not |[N_{G}(P):P]$. In the second case, $N_{G_{p}}(P)$ is a quaternion group of order eight 
and $4\not |[N_{G}(P):P]$. If $P$ is generalized quaternion and not a Sylow $p$-subgroup, 
it follows that $N_{G_{p}}(P)$ is a generalized quaternion group with 
$[N_{G_{p}}(P):P]=2$ while if $P$ is a Sylow $p$-subgroup, $[N_{G_{p}}(P):P]=1$. 
We have now proven in each case that if $P$ is a $p$-centric subgroup of 
$G$, then $p^2\not |[N_{G}(P):P]$; furthermore, the map $\alpha _{p}$ is a bijection. 
From our previous discussion, we see that the map $\phi _{G,p}$ is also a bijection.
\end{proof}

The proof of Theorem \ref{thm:odd} follows from the next lemma.
\begin{lemma}[See {\cite[6.3.10]{s:gtt}}]\label{lem:odd}
Let $G$ be a finite group of odd order. If $p||G|$ and $P$ is a 
principal $p$-radical subgroup of $G$ such that $\rk(P)\leq 2$, then $P$ 
is a Sylow $p$-subgroup of $G$.
\end{lemma}
\begin{proof} Let $p||G|$ and $P$ be a principal $p$-radical subgroup of $G$ such that $\rk(P)~\leq~2$. 
Suppose $P$ is not a Sylow $p$-subgroup of $G$. Let $W=N_{G}(P)/C_{G}(P)$. The derived 
subgroup $W'$ of $W$ is a $p$-subgroup (see {\cite[6.3.10]{s:gtt}}). For a subgroup 
$H$ of $N_{G}(P)$, we will denote by $\bar{H}$ the image of $HC_{G}(P)$ in $W$. Now 
since $P$ is principal $p$-radical, $O_{p}(W)=\bar{P}$. This implies $W' \subseteq 
\bar{P}$. So $W/\bar{P}=N_{G}(P)/PC_{G}(P)$ is an abelian group. Since $P$ is not 
a Sylow $p$-subgroup of $G$, $p|[N_{G}(P):P]$. $P$ is principal $p$-radical; therefore, 
$Z(P)$ is a Sylow $p$-subgroup of $C_{G}(P)$. This shows that $p|[N_{G}(P):PC_{G}(P)]$; 
thus, $N_{G}(P)/PC_{G}(P)$ has a non-trivial normal $p$-subgroup. This contradicts 
$P$ being a principal $p$-radical subgroup of $G$, so $P$ is a Sylow $p$-subgroup of 
$G$. \end{proof}

According to Lemma \ref{lem:odd}, if $G$ is a rank two group of odd order
with $p||G|$, then the map $\alpha _{p}$ 
is a bijection, and the map $\phi _{G,p}$ is also a bijection. It follows from the 
introduction that if $G$ is a rank two group of odd order, then the map $\phi _{G}$ is also a 
bijection, thus proving Theorem \ref{thm:odd}.\\
\\
This only leaves Theorem \ref{thm:rk2} to be proven. As we have been doing in order to 
show that the map $\phi _{G}$ is a surjection, we will be showing that the map $\phi _{G,p}$ is a 
surjection for 
each prime $p$ dividing $|G|$. We will need to approach the prime $2$ differently then the odd 
primes. Thus, we will break Theorem \ref{thm:rk2} into the following two propositions.
\begin{proposition}\label{prop:odd}\fussy
Let $G$ be a finite group and $p>2$ an odd prime dividing $|G|$. If $G$ does not contain a rank 
three elementary abelian $p$-subgroup, then $\phi _{G,p}: [BG,BU(n)^{\wedge}_{p}] \rightarrow
\Char _{n}^{G}(G_{p})$ is a surjection.
\end{proposition}

\begin{proposition}\label{prop:p2}
If $G$ be a finite group of even order not containing a rank three elementary abelian 
$2$-subgroup, then $\phi _{G,2}: [BG,BU(n)^{\wedge}_{2}] \rightarrow \Char _{n}^{G}(G_{2})$ 
is a surjection.
\end{proposition}
It is clear that 
Theorem \ref{thm:rk2} follows immediately from these two propositions. We will spend the 
next two sections proving these two propositions.
\section{$3$-Groups of Maximal Class}\label{sec:max}
We will prove Proposition \ref{prop:odd} first. Before we do, we 
need to focus on a type of group that 
will arise in this proof. For Section \ref{sec:max} we will focus on the case where a 
Sylow $3$-subgroup of $G$ is a $3$-group of maximal class with the purpose of
proving the following lemma. 
\begin{lemma}
Let $G$ be a finite group and $G_{3}$ be a Sylow $3$-subgroup of $G$ such that 
$G_{3}$ is a $3$-group of maximal class and $|G_{3}|\geq 3^5$. 
If $P\subseteq G_{3}$ is a $3$-centric subgroup of 
$G$, then $[N_{G_{3}}(P):P]\leq 3^2$.
\end{lemma}
\begin{proof} We will let the series 
$$Z(G_{3}) = C_{n-1}\subset C_{n}\subset \cdots \subset C_{3}\subset
C_{2}\subset G_{3}$$ be the lower central series of $G_{3}$. This is defined inductively 
by $C_{2}=[G_{3},G_{3}]$ and $C_{i}=[C_{i-1},G_{3}]$ for $i>2$. Since $G_{3}$ is of 
maximal class, $n$ is the integer such that $|G_{3}|=3^n$. We also 
introduce another subgroup of $G_{p}$, which will be denoted by $C_{1}$.  
The subgroup $C_{1}$ is defined by the property that $C_{1}/C_{4}$ is 
the centralizer in $G_{3}/C_{4}$ of $C_{2}/C_{4}$ (see {\cite[section 2]{b:scpg}}). 
By these definitions we note that each $C_{i}$ is characteristic in $G_{p}$, thus giving the 
increasing sequence of $n-1$ distinct proper subgroups of $G_{p}$:
$$Z(G_{3}) = C_{n-1}\subset C_{n}\subset \cdots \subset C_{3}\subset
C_{2}\subset C_{1}\subset G_{3}.$$
By the work of Blackburn \cite{b:scpg}, either $C_{1}$ is abelian or 
the commutator subgroup $C_{2}$ is abelian and $[C_{1},C_{2}]=Z(G_{3})$. Now assume that 
$P$ is a proper subgroup of $G_{3}$ that is a $3$-centric subgroup of $G$. We will look at 
two cases: the first when $P\subseteq C_{1}$ and the second when $P\not \subseteq C_{1}$.\\
\\
Case I: Assume that $P\subseteq C_{1}$. Since $P$ is a $3$-centric subgroup of $G$, 
$Z(C_{1})$ must be a proper subgroup of $P$. Notice that if either $C_{1}$ is abelian or 
the commutator subgroup $C_{2}$ is abelian and $[C_{1},C_{2}]=Z(G_{3})$, then 
$C_{3}\subseteq Z(C_{1})$. This inclusion implies that $[G_{3}:P]\leq 3^2$.\\
\\
Case II: Assume that $P\not \subseteq C_{1}$. Let $N=N_{G_{3}}(P)$. Let $t$ be 
an element of $P\setminus (P\cap C_{1})$. It can be shown that $|C_{G_{3}}(t)|=3^2$ 
and that $C_{G_{3}}(t)\subseteq P$; therefore, any subgroup of $G_{3}$, which contains 
$P$, is of maximal class.\\
\\
Suppose that $P$ is abelian. This implies that $P=C_{G_{3}}(t)$ and  that 
$N$ is generated by $t$ and $C_{n-2}$; thus $[N:P]=3$.\\
\\
Now suppose that $P$ is not abelian. Since $N$ is of maximal class, 
$[N,N]\subseteq C_{2}$; therefore, $[N,N]$ is abelian. Seeing that $P\lhd N$ and $P$ is not 
abelian, $[N:P]=3$.\end{proof}

\section{Finite groups, part II}

We need to mention one more lemma before we proceed with the proof of Proposition 
\ref{prop:odd}. This Lemma deals with the case where a Sylow $p$-subgroup of $G$ is a 
metacyclic $p$-group and $p$ is an odd prime. It follows from the work of J. Dietz \cite{d:sscs} 
and of Martino and Priddy \cite{mp:chsg}.
\begin{lemma}[See \cite{d:sscs}, and {\cite[Theorem 2.7]{mp:chsg}}]\label{lem:meta}
Let $G$ be a finite group, $p$ an odd prime divisor of $|G|$, and $G_{p}$ a Sylow 
$p$-subgroup of $G$. If $G_{p}$ is a metacyclic group, $|G_{p}|>p^3$, and $P \subseteq G_{p}$ 
is a principal $p$-radical subgroup of $G$, then $P=G_{p}$.
\end{lemma}
\begin{proof} Suppose $G_{p}$ is a metacyclic group with $|G_{p}|> p^{3}$. Martino and Priddy 
\cite{mp:chsg} have shown that such a $G_{p}$ is a Swan Group. In the proof they proved that if 
$H\subset G_{p}$ is a proper subgroup with $C_{G_{p}}(H) =Z(H)$, then 
\[(N_{G_{p}}(H)/HC_{G_{p}}(H))\cap O_{p}(\out(H))\neq \{1\}.\] 
So no proper subgroup $H\subset G_{p}$ can be a principal $p$-radical subgroup of $G$.
\end{proof}

Now we give the proof of Proposition \ref{prop:odd}:
\begin{proof} Recall that $G$ is a finite group that does not contain a rank three elementary abelian 
subgroup and $p$ is an odd prime dividing $|G|$. Let $G_{p}$ be a Sylow $p$-subgroup 
of $G$ and let $n$ be an integer such that $|G_{p}|=p^n$. Since Theorem \ref{thm:rk1} takes care 
of the case where $\rk (G_{p})=1$, we may assume that $\rk (G_{p})=2$. It has been shown by 
Blackburn \cite{b:gcet} that one of the following holds for $G_{p}$ (see also \cite{dp:shtr2}):
\begin{enumerate}
\item $n<5$,
\item $n\geq 5$ and $G_{p}$ is metacylic,
\item $n\geq 5$, $p=3$, and $G_{3}$ is a $3$-group of maximal class,
\item $n\geq 5$ and $G_{p}= \langle a,b,c|a^p=b^p=c^{p^{n-2}}=1, [a,b]=c^{p^{n-3}}, 
c \in Z(G_{p})\rangle$, or
\item $n\geq 5$, $e\neq 0$ is a quadratic nonresidue mod $p$, 
and \newline $G_{p}= \langle a,b,c|a^p=b^p=c^{p^{n-2}}=[b,c]=1, [a,b^{-1}]=c^{e p^{n-3}}, 
[a,c]=b\rangle$.
\end{enumerate}
Recall that it is enough to show that in each of these cases, if $P\subseteq G_{p}$ is a principal
$p$-radical subgroup of $G$, then $[N_{G_{p}}(P):P]\leq p^2$. Lemma \ref{lem:meta} gives 
this result for case 2, while case 3 was shown in Section \ref{sec:max}. To proceed with the other 
three cases, recall that by Lemma \ref{lem:cent} it is enough to show that 
$[G_{p},Z(G_{p})]\leq p^3$. The first case follows immediately since $G_{p}$ must have a 
non-trivial center. Notice in both of the remaining cases that the center of $G_{p}$ contains the 
element $c^{p}$ implying that $|Z(G_{p})|\geq p^{n-3}$. 
This concludes these two cases as well as the proof of Proposition \ref{prop:odd}.\end{proof}
We are also ready to give the proof of of Proposition \ref{prop:p2}:
\begin{proof} Recall that $G$ is a finite group of even order that does not contain a rank three elementary 
abelian subgroup. As before let $G_{2}$ be a Sylow $2$-subgroup of $G$. Suppose that 
$P\subset G_{2}$ is a principal $2$-radical subgroup of $G$, which is not itself a Sylow 
$2$-subgroup of $G$. Notice that $N_{G}(P)/P C_{G}(P)\subseteq \out(P)$ is both non-trivial 
and not a $2$-group. In particular this means that $P$ is a $2$-group of rank one or two 
with an automorphism of odd order. Richard Thomas has classified all such $2$-groups 
\cite{t:apo5,t:apo3}. According to his classification $P$ must be one of the following:
\begin{enumerate}
\item $P\cong Q_{8}\ast D_{8}$,
\item $P\cong Q_{8}\times Q_{8}$,
\item $P\cong Q_{8}\wr \mathbb{Z}_{2}$,
\item $P\cong Q_{8}\ast C$ where $\out (C)$ is a $2$-group, $\out(C)\lhd \out(P)$, 
and $C\neq D_{8}$,
\item $P\cong \mathbb{Z}_{2^{r}} \times \mathbb{Z}_{2^{r}}$, or
\item $P\cong U_{64}\in \Syl _{2}(\mathrm{PSU}_{3}(\mathbb{F}_{4}))$.
\end{enumerate}
Our notation is as follows: $Q_{8}$ means the quaternion 8 group, $D_{8}$ means the 
dihedral group of size 8, $\ast$ refers to the central product, and $\wr$ refers to the wreath 
product.\\
\\
Since $P$ is a principal $2$-radical subgroup of $G$, we know from the definition that
$N_{G}(P)/P C_{G}(P)\cap O_{2}(\out(P))=\{1\}$; 
therefore, let $S\in \Syl _{2}(N_{G}(P)/P C_{G}(P))$ and $T\in \Syl _{2}(\out(P)/O_{2}(\out(P)))$.
It is clear then that $S$ must be isomorphic to a subgroup of $T$. In order to show that 
$\rk _{2}(N_{G}(P)/P)\leq 2$, it is enough to show that $\rk _{2}(T)\leq 2$. We will do this on a 
case by case basis using the list above. In cases 1 and 2, $T\cong D_{8}$, giving 
$\rk _{2}(T)=2$. In cases 3, 4, and 5, $T\cong \mathbb{Z}_{2}$, implying that $\rk _{2}(T)=1$. 
In the final case, $T\cong \mathbb{Z}_{4}$, implying that $\rk _{2}(T)=1$. We see that for any 
possible $P\subseteq G_{2}$ that is a principal $2$-radical subgroup of $G$, 
$\rk _{2}(N_{G}(P)/P)\leq 2$. 
This concludes the proof of Proposition \ref{prop:p2}. \end{proof}

\end{document}